\documentstyle[12pt,amssymb,amstex]{amsart}

\numberwithin{equation}{section}

\def\a{\alpha}

\def\g{\gamma}

\def\e{\epsilon}
\def\l{\lambda}

\def\s{\sigma} 
\def\t{\tau}

  \def\Th{\Theta}

\newtheorem{thm}{Theorem}[section]
\newtheorem{lem}[thm]{Lemma}

\begin{document}

\small

\title[The almost everywhere convergence of the spectral
expansions] {About the almost everywhere convergence of the
spectral expansions of functions from $L_1^\a(S^N)$ }

\author{Anvarjon Akhmedov}
\email{{\tt anv.akhmed@@gmail.com},{\tt
anvarjon\_ahmedov@@yahoo.com}}

\begin{abstract}
In this paper we study  the almost everywhere convergence of the
expansions related to the self-adjoint extension of the
Laplace-Beltrami operator on the unit sphere. The sufficient
conditions for summability is obtained. The more general
properties and representation by the eigenfunctions of the
Laplace-Beltrami operator of the Liouville space $L_1^\a$ is used.
 For the orders of Riesz means, which greater than
critical index $\frac{N-1}{2}$ we proved the positive results on
summability of Fourier-Laplace series.  Note that when order
$\alpha$ of Riesz means is less than critical index then for
establish of the almost everywhere convergence requests to use
other methods form proving negative results. We have constructed
different method of summability of Laplace series, which based on
spectral expansions property of self-adjoint Laplace-Beltrami
operator on the unit sphere.
\end{abstract}

\maketitle

\footnotetext[1]{The author (A.A.) is on leave from National
University of Uzbekistan, Department of Mechanics \& Mathematics,
Vuzgorodok, 100174 Tashkent,  Uzbekistan }

\section{Introduction}

Let $S^N$ is unit sphere in $R^N$
$$
S^N=\{x\in R^N: |x|=1\}
$$
The sphere $S^N$ is naturally equipped with a positive measure
$d\s(x)$ and with an elliptic second order differential operator
$\Delta_s$, which named the Laplace-Beltrami operator on the
sphere. This operator is symmetric and nonnegative, it extends to
a nonnegative self-adjoint operator on the space $L_2(S^N)$ (by
$L_p(S^N)$ we mean the $L_p$-space associated with the measure
$d\s(x)$ on the sphere). Let $c$ be any positive number, and let
$A=\Delta_s+1$. We denote by $Spec(A)=\{\l_k, k=0,1,2,...\}$ the
spectrum of $A$. This spectrum is nondecreasing sequence of
positive eigenvalues with finite multiplicities (and written as
such) tending to infinity. We denote by $Y_j^k(x)$ an
eigenfunction of the Laplace-Beltrami Operator.
$$\Delta_sY_j^{(k)}=\l_kY_j^{(k)}$$
where $\l_k=k(k+N-1),k=0,1,2....$ The system of eigendunctions of
the Laplace-Beltrami operator is orthonormal basis in
$L_2^(S^N).$\\
 One of the main problems of harmonic analysis is
the reconstruction of functions from their expansion:
\begin{equation}\label{1}
f(x)\sim\sum\limits_{k=0}^\infty Y_k(f,x)
\end{equation}
 The main purpose this article the
convergence problems of the partial sums of representation \ref{1}
\begin{equation}\label{eq1}
E_nf(x)=\sum\limits_{k=0}^nY_kf(x,y)
\end{equation}
in (\ref{1}) and (\ref{eq1}) denoted by $Y_k(f,x)$
$$Y_k(f,x)=\int\limits_{S^N}f(y)Z_k(x,y)d\sigma(y)$$
where $Z_k(x,y)=\sum\limits_{k=o}^nY_j^{(k)}(x)Y_j^{(k)}(y)$ is
Zonal harmonic of order $k.$ A spectral function $\Th(x,y,\l)$  of
the Laplace-Beltrami operator on sphere is defined by
\begin{equation}\label{eq2}
\Th(x,y,n)=\sum\limits_{k=0}^nZ_k(x,y)
\end{equation}
Then the spectral expansions of the (\ref{eq1}) can be rewritten
as follow
\begin{equation}
E_nf(x)=\int\limits_{S^N}f(y)\Th(x,y,n)d\s
\end{equation}
The Riesz means of the spectral expansions (\ref{eq1}) is defined
buy the next expression
\begin{equation}\label{eq3}
E_n^\a f(x)=\int\limits_{S^N}f(y)\Th^\a(x,y,n)d\s
\end{equation}
where by $\Th^\a(x,y,\l)$ denoted the Riesz means of the spectral
function (\ref{eq2}):
\begin{equation}
\Th^\a(x,y,n)=\sum\limits_{k=0}^n\left(1-\frac{\l_k}{\l_n}\right)^\a
Z_k(x,y)
\end{equation}
Let us denote by $T_n^\a f(x)$ the Cesaro means of the spectral
expansions, which defined as follow
\begin{equation}
T_n^\a f(x)=\int\limits_{S^N}f(y)\Phi^\a(x,y,n)d\s
\end{equation}
with the kernel
\begin{equation}
\Phi^\a(x,y,n)=\sum\limits_{k=0}^n\frac{A_{n-k}^\a}{A_n^\a}Z_k(x,y)
\end{equation}
where
$$A_m^\a=\frac{(\a+1)(\a+2)...(\a+m)}{m!}
$$
The maximal operator of the Integral operator (\ref{eq3}) plays
important role in establishing the almost everywhere convergent of
spectral expansions
$$E_\star^\a f(x)=\limsup_{n\in N}|E_n^\a f(x)|$$
Let $L_p^\t(S^N)$ is the class Lioville, which consists of all
functions from $L_p(S^N)$, for which
\begin{equation}
\|f\|_{L_p^\t}=\|\sum\limits_{k=0}^\infty\l_k^\t
Y_k(f,x)\|_{L_p}<\infty
\end{equation}
Main result of this paper is

\begin{thm}\label{theo}
Let $f\in L_1^\t(S^N),
\a+\t>\frac{N-1}{2},0\leq\a\leq\frac{N-1}{2},\t>0,$ then\\
1) $\|E_{\star}^\a f(x)\|_{L_1(S^N)}\leq
c_\a\|f\|_{L_1^\t(S^N)}$\\
2) The Riesz means $E_n^\a f(x)$ almost everywhere on $S^N$
converges to  $f(x).$
\end{thm}

Note that, if $\t>N$ then Riesz means of spectral expansions
uniformly convergence to $f$ (see \cite{Al}). In work \cite{Bas}
was investigated spectral expansions related to the
Pseudo-differential operators and in the \cite{Buv} considered
spectral expansions of elliptic differential operators.

\section{Preliminaries}

In this section we recall some preliminaries concerning some
properties of Liuoville space and Riesz means.

\begin{lem}\label{emb}
Let $\t\geq 0,\a\geq 0, p\geq 1,$ then for all $f\in
C^\infty(S^N)$ we have
 $$\|A^{\a/2}f\|_{L_p(S^N)}\leq c_\a\|f\|_{L_p^{\a+\t}(S^N)}$$
\end{lem}

The proof of this lemma it follows from embedding properties of
Liouvill space. For more details see \cite{Al}.

\begin{lem}\label{sp}
Let $\Th^\a(x,y,n)$ is the Riesz means of the spectral function of
the Laplace-Beltrami operator on sphere and\\
 1) if spherical
distance $\g=\g(x,y)$ between $x$ and $y$ satisfied inequality
$|\frac{\pi}{2}-\g|<\frac{n}{n+1}\frac{\pi}{2}$ then we have
$$
\Th^\a(x,y,n)=O(1)(\frac{n^{(N-1)/2}}{(\sin\g)^{(N-1)/2}(\sin(\g/2))^{1+\a}}+
\frac{n^{(N-3)/2}}{(\sin\g)^{(N+1)/2}(\sin(\g/2))^{1+\a}}+$$
\begin{equation}+\frac{n^{-1}}{(\sin(\g)/2)^{1+N}});
\end{equation}
2) if $0\leq\g\leq\pi$, then we have $$\Th^\a(x,y,n)=O(1)n^N;$$

 3) if $0<\g_0\leq\g\leq\pi$, then we have $$\Th^\a(x,y,n)=O(1)n^{N-\a};$$
\end{lem}
This lemma proved in our work \cite{Akh}.

\begin{lem}\label{Rz}
Let $f\in C^\infty(S^N),\a+\t>\frac{N-1}{2}
0\leq\a\leq\frac{N-1}{2},\t>0$ then for maximal operator of Riesz
means we have
 $$\|E_{\star}^\a f\|_{L_1(S^N)}\leq c_\a\|f\|_{L_1^\t(S^N)}$$
\end{lem}
Proof. For all $\t\geq 0$ let us denote by $\Th_\t^\a(x,y,n)$ the
kernel of the integral operator $A^{-\t}E_n^\a$:
$$\Th_\t^\a(x,y,n)=\sum\limits_{k=1}^n\l_k^{-\t}\left(1-\frac{\l_k}{\l_n}\right)^\a
Z_k(x,y)$$ Let $g=A^{\t/2}f$. Using the equality
$$E_n^\a f=A^{-\t/2}E_n^\a A^{\t/2}f$$we can rewrite the Riesz
means as follow
$$
E_n^\a f(x)=\int\limits_{S^N}\Th_{\t/2}^\a(x,y,n)g(y)d\s(y). $$

The kernel $\Th_{\t/2}^\a(x,y,n)$
 we reduce to the next form
$$\Th_{\t/2}^\a(x,y,n)=\sum\limits_{k=1}^n\l_k^{-\t/2}
\left(\Th^\a(x,y,k)-\Th^\a(x,y,k-1)\right)=$$
$$=\sum\limits_{k=1}^n\left(\l_k^{-\t/2}-\l_{k+1}^{-\t/2}\right)
\Th^\a(x,y,k)+\l_n^{-\t/2}\Th^\a(x,y,n).$$

Using this representation for the kernel $\Th_{\t/2}^\a(x,y,n)$ we
can estimate the Riesz means of the spectral expansions as follow:
$$E_n^\a f(x)=\sum\limits_{k=1}^n\left(\l_k^{-\t/2}-\l_{k+1}^{-\t/2}\right)
\int\limits_{S^N}\Th^\a(x,y,k+1)g(y)d\s(y)+$$
$$+\l_n^{-\t/2}\int\limits_{S^N}\Th^\a(x,y,n)g(y)d\s(y).$$
Separating into four part the integral
$\int\limits_{S^N}\Th^\a(x,y,k+1)g(y)d\s(y)$ and for estimate each
part let us apply the lemma \ref{sp}:
$$\int\limits_{S^N}\Th^\a(x,y,k)g(y)d\s(y)=
I_1+I_2+I_3+I_4,$$ where

$$I_1=\int\limits_{\g<\frac{1}{k}}\Th^\a(x,y,k)g(y)d\s(y)$$
$$I_2=\int\limits_{\frac{1}{k}<\g\leq\frac{\pi}{2}}\Th^\a(x,y,k)g(y)d\s(y)$$
$$I_3=\int\limits_{\frac{\pi}{2}<\g\leq\pi-frac{1}{k}}\Th^\a(x,y,k)g(y)d\s(y)$$
$$I_4=\int\limits_{\pi-\frac{1}{k}<\g\leq\pi}\Th^\a(x,y,k)g(y)d\s(y).$$

For estimate $I_1$ and $I_4$ we use the asymptotic behavior
\ref{sp} of the Riesz means of the spectral function:
$\Th^\a(x,y,n)=O(1)n^N$ and we get
$$|I_1+I_4|\leq C\left(k^N\int\limits_{\g(x,y)<\frac{1}{k}}|g(y)|d\s(y)+
k^N\int\limits_{\g(\overline{x},y)=\pi-\g(x,y)<\frac{1}{k}}|g(y)|d\s(y)\right)\leq$$
$$\leq\left(g^{\star}(x)+g^{\star}(\overline{x})\right)$$
here we denoted by $g^{\star}(x)$ the maximal function of the
function $d(x)$ which defined by the formula
$$g^{\star}(x)=\limsup_{r>}\frac{1}{mesB(x,r)}\int\limits_{B(x,r)}|g(y)|d\s(y)$$
where $B(x,r)$ is the ball in the unit sphere with the center at
the point $x$ and the radius $r$:
$$B(x,r)=\{y\in S^N: \g(x,y)<r, r>0\}.$$
Due to estimation 1) of lemma \ref{sp} we have for $I_2$
$$|I_2|\leq
Ck^{\frac{N-1}{2}-\a}\int\limits_{\frac{1}{k}<\g\leq\frac{\pi}{2}}(\sin\g)^{(-\frac{N+1}{2}-\a}|g(y)|d\s(y)+$$
$$+Ck^{\frac{N-3}{2}-\a}\int\limits_{\frac{1}{k}<\g\leq\frac{\pi}{2}}(\sin\g)^{(-\frac{N+3}{2}-\a}|g(y)|d\s(y)+$$
$$+Ck^{-1}\int\limits_{\frac{1}{k}<\g\leq\frac{\pi}{2}}(\sin\g)^{-1-N}|g(y)|d\s(y)\leq
Cg^{\star}(x)\left(1+k^{\frac{N-1}{2}-\a}\right)$$

The $I_3$ after the denoting $\g(\overline{x},y)=\pi-\g(x,y)$ can
be computed as $I_2$, and we have
$$|I_3|\leq Cg^{\star}(\overline{x})\left(1+k^{\frac{N-1}{2}-\a}\right)$$
Such that for Riesz means of the spectral expansions we have
$$E_n^\a f(x)\leq C\sum\limits_{k=1}^n\left(\l_k^{-\t/2}-\l_{k+1}^{-\t/2}\right)
\left(1+k^{\frac{N-1}{2}-\a}\right)\left(g^{\star}(x)+g^{\star}(\overline{x})\right)+$$
$$+\l_n^{-\t/2}\left(1+n^{\frac{N-1}{2}-\a}\right)\left(g^{\star}(x)+g^{\star}(\overline{x})\right).$$
We have to prove that the expression
$$T_n=\sum\limits_{k=1}^n\left(\l_k^{-\t/2}-\l_{k+1}^{-\t/2}\right)
\left(1+k^{\frac{N-1}{2}-\a}\right)
$$
is bounded with constant which does not depend of $n.$ If note
that eigenvalues $\l_k=k(k+N-1), k=1,2,3,...$

$$T_n=\sum\limits_{k=1}^n\left(k^{-\t/2}(k+N-1)^{-\t/2}-(k+1)^{-\t/2}(k+N)^{-\t/2}\right)
k^{\frac{N-1}{2}-\a}=
$$

$$
=\sum\limits_{k=1}^n\left(\left(1+\frac{N-1}{k}\right)^{-\t/2}-\left(1+\frac{N+1}{k}+\frac{N}{k^2}\right)^{-\t/2}\right)
k^{\frac{N-1}{2}-\a-\t}
$$

and applying the Lagrange formula to function $J(x)=(1+x)^{-\t/2}$
in the segment $[\frac{N-1}{k},\frac{N+1}{k}+\frac{N}{k^2}]$ then
we obtain
$$J\left(\frac{N-1}{k}\right)-J\left(\frac{N+1}{k}+\frac{N}{k^2}\right)=
-\left(\frac{2}{k}+\frac{N}{k^2}\right)J'(\xi)
$$
where $\xi: \frac{N-1}{k}<\xi<\frac{N+1}{k}+\frac{N}{k^2}$. Then
for Riesz means we have

\begin{equation} \label{c}
|E_n^\a f(x)|\leq C =\sum\limits_{k=1}^\infty
k^{\frac{N-1}{2}-\a-\t-1}
\end{equation}

If $\a+\t>\frac{N-1}{2}$ then numerical series \ref{c} converges.
Such that for Riesz means of order $\a$ we have proved estimation
\begin{equation}\label{est}
|E_{\star}^\a f(x)|\leq C
\leq\left(g^{\star}(x)+g^{\star}(\overline{x})\right)
\end{equation}
for all $f\in L_1^\t(S^N)$ provided by the condition
$\a+\t>\frac{N-1}{2}.$

Let $p>1,$ then
$$\|E_{\star}^\a f\|_{L_1}\leq C\|E_{\star}^\a f\|_{L_p}\leq\ C
\left(\|g^{\star}(x)\|_{L_p}+\|g^{\star}(\overline{x})\|_{L_p}\right)\leq
C\|g\|_{L_p}
$$
By virtue of the lemma \ref{emb} we have
$$\|g\|_{L_p}\leq C\|f\|_{L_p^\a}$$
It is well-known the embedding $L_{p+N(1-\frac{1}{p})+\e}\subset
L_1^\a$, if $\e>0$ and $p>1.$ Therefore finally we get the
estimate \begin{equation} \label{eq5} \|E_{\star}^\a
f\|_{L_1(S^N)}\leq c_\a\|f\|_{L_1^\t(S^N)}
\end{equation}

 The assumption of lemma \ref{Rz} is proved for all function $f\in
C^\infty$. It is well-known, that Banach space $C^\infty$ is dense
in $L_1^\t$ for all $\t>0.$ Therefore  proof of the inequality
\ref{eq5} for all function from $L_1^\t$ follows from density of
$C^\infty$ in $L_1^\t.$ Lemma \ref{Rz} is proved.

Let us prove second part of the theorem\ref{theo}. Let $f\in
L_1^\t(S^N),\t>0,\a+\t>\frac{N-1}{2}.$ Fix any $\e>0.$ From
density of $C^\infty$, it follows that exists $h\in C^\infty$ such
that $\|f-h\|_{L_1^\t(S^N)}<\e.$ Riesz means $E_n^\a h(x)$ of
spectral expansions of $h$ uniformly converges to $h(x).$ For any
$\e>0$ there is $n_0(\e)$ such that, if $n>n_0$ we have
$$|E_n^\a h(x)-h(x)|<\e, x\in S^N.$$
Then for $n>n_0$ we have
$$|E_n^\a f(x)-f(x)|<E_{\star}^\a (f-h)(x)+\e+|h(x)-f(x)|.$$
Using first part of the theorem \ref{theo} it is not hard to see
that
$$\|E_{\star}^\a f(x)-f(x)\|_{L_1}\leq \|h(x)-f(x)\|_{L_1^\t(S^N)}+\e.$$
The last inequality proves that
$$ \lim_{n\rightarrow\infty}E_n^\a f(x)=f(x)$$
holds almost everywhere in $S^N.$

For Cesaro means the assertion of the theorem \ref{theo} proves
analogously.


\begin{thebibliography}{15}%EIBLM}

\bibitem{Sh}
M. A. Shubin, Pseudodifferential operators and spectral Theory,
{\it Moscow. Science Publisher} 1978, pp.280.

\bibitem{Al}
Sh.A. Alimov,   {\it Journal of Differential Equations.} T. IX,
V.4, 669-681,(1973).

\bibitem{Bas}
A.Y. Bastis,  {\it Mathematical Notes}, T.34, V.4,587-600,(1983).

\bibitem{Buv}
K.T. Buvaev,   {\it Uzbek Mathematics Journal}V.3, 39-44,(1995).

\bibitem{Kog}
J. Kogbetlians , {\it Journal of Mathematical Pure et Apple} ,Ser
9,1924

\bibitem{Akh}
A.A. Akhmedov ,   {\it Uzbek Mathematics Journal}V.4,
17-26,(1996).

%\bibliographystyle{amsplain}
\end{thebibliography}
\end{document}